\documentclass[11pt,reqno]{amsart}
\usepackage{amsfonts}
\usepackage{mathrsfs}
\usepackage{amsmath}
\usepackage{amssymb,amsfonts}

\newtheorem{thm}{Theorem}[section]

\newtheorem{coro}{Corollary}[section]\numberwithin{equation}{section}
\newtheorem{rmk}{Remark}[section]

\def\pf{{\textit {Proof:} }}

\newcommand{\mysection}[1]{\section{#1}\setcounter{equation}{0}}

\newfont{\bb}{msbm10 at 11pt}

\newcommand{\bal}{\begin{aligned}}      \newcommand{\eal}{\end{aligned}}
\newcommand{\ba}{\begin{array}}      \newcommand{\ea}{\end{array}}
\newcommand{\bc}{\begin{center}}     \newcommand{\ec}{\end{center}}
\newcommand{\be}{\begin{enumerate}}  \newcommand{\ee}{\end{enumerate}}
\newcommand{\beq}{\begin{eqnarray}}  \newcommand{\eeq}{\end{eqnarray}}
\newcommand{\beQ}{\begin{eqnarray*}} \newcommand{\eeQ}{\end{eqnarray*}}
\newcommand{\bi}{\begin{itemize}}    \newcommand{\ei}{\end{itemize}}
\newcommand{\bt}{\begin{tabular}}    \newcommand{\et}{\end{tabular}}
\newcommand{\bdm}{\begin{displaymath}} \newcommand{\edm}{\end{displaymath}}

\def\qed{\hfill{Q.E.D.}\smallskip}

\newcommand{\ls}{\setlength{\baselineskip}{12pt}
                 \setlength{\parskip}{3mm}}

\begin{document}

\title[Eguchi-Hanson Types]{Metrics of Eguchi-Hanson types with the negative constant scalar curvature}

\author[J Chen]{Junwen Chen$^{\dag}$}
\address[]{$^{\dag}$College of Mathematics and Information Science, Guangxi University, Nanning, Guangxi 530004, PR China}
\email{jwchen20@st.gxu.edu.cn}
\author[X Zhang]{Xiao Zhang$^{\flat}$}
\address[]{$^{\flat}$Guangxi Center of Mathematical Research, Guangxi University, Nanning, Guangxi 530004, PR China}
\address[]{$^{\flat}$Institute of Mathematics, Academy of Mathematics and Systems Science, Chinese Academy of Sciences,
Beijing 100190, PR China and School of Mathematical Sciences, University of Chinese Academy of Sciences, Beijing 100049, PR China}
\email{xzhang@gxu.edu.cn, xzhang@amss.ac.cn}

\date{}

\begin{abstract}
We construct two types of Eguchi-Hanson metrics with the negative constant scalar curvature. The type I metrics are K\"{a}hler. The type II metrics are ALH whose total energy can be negative.
\end{abstract}

\maketitle \pagenumbering{arabic}

\mysection{Introduction}\ls

Eguchi-Hanson metrics are Ricci flat, anti-self-dual 4-dimensional ALE Riemannian metrics \cite{EH1,EH2}. They are referred as gravitational instantons, and attract much attention in general relativity. The counter-examples of Hawking and Pope's generalized positive action conjecture, which were first constructed by LeBrun using the method of algebraic geometry, are Burns' metrics and have an opposite orientation of Eguchi-hanson metrics \cite{L1}. The direct construction of metrics of Eguchi-Hanson type with the zero scalar curvature, negative mass was provided in \cite{Z} by solving an ordinary differential equation.

The Eguchi-Hanson metrics with constant Ricci curvatures and their relations with K\"{a}hler geometry were studied by Pedersen \cite{P}. The general condition for metrics of Eguchi-Hanson type to be K\"{a}hler was derived by LeBrun \cite{L2}. In this short paper, we construct two types of Eguchi-Hanson metrics with the negative scalar curvature. The metrics of type I have natural K\"{a}hler structure, but not ALH. They can not develop into 5-dimensional spherically symmetric, static vacuum Einstein spacetimes with the negative cosmological constant. The metrics of type II are not K\"{a}hler, but ALH which allow negative total energy. They can only develop into 5-dimensional spherically symmetric, static vacuum Einstein spacetimes with the negative cosmological constant as Clarkson-Mann's metrics. These metrics are interesting in the theory of AdS/CFT correspondence.

The paper is organized as follows. In Section 2, we provide a brief introduction to metrics of Eguchi-Hanson type with the zero scalar curvature. In Section 3, we construct metrics of Eguchi-Hanson type I with the negative constant scalar curvature. In Section 4, we construct metrics of Eguchi-Hanson type II with the negative constant scalar curvature. In Section 5, we compute the total energy of metrics of Eguchi-Hanson type II and show that it can be negative.

\mysection{Eguchi-Hanson Metrics}\ls

Let $\theta$, $\phi$, $\psi$ be the Euler angles on the $3$-sphere $S^3$ with ranges
$$0 \leqslant \theta \leqslant \pi, \quad 0 \leqslant \phi \leqslant 2\pi, \quad 0 \leqslant \psi \leqslant 4\pi. $$
Let $\sigma_{1}$, $\sigma_{2}$, and $\sigma_{3}$ be the Cartan-Maurer one-forms for $SU(2)\cong S^3$, given by
\begin{equation*}
\begin{aligned}
\sigma_{1}&=\frac{1}{2}\big(\sin\psi d\theta-\sin\theta cos\psi d\phi \big),\\
\sigma_{2}&=\frac{1}{2}\big(-\cos\psi d\theta-\sin\theta \sin\psi d\phi \big),\\
\sigma_{3}&=\frac{1}{2}\big(d \psi+\cos\theta d\phi \big).
\end{aligned}
\end{equation*}
They satisfy $$d\sigma_{1}=2\sigma_{2}\wedge\sigma_{3},\quad d\sigma_{2}=2\sigma_{3}\wedge  \sigma_{1},\quad d\sigma_{3}=2\sigma_{1}\wedge  \sigma_{2}.$$

In \cite{EH1, EH2}, Eguchi and Hanson constructed 4-dimensional self-dual Riemannian metrics
 \beq
g =\Big(1-\frac{B}{r ^4}\Big)^{-1} dr ^2+r ^2 \Big( \sigma _1
^2 +\sigma _2 ^2 +\big(1-\frac{B}{r ^4} \big)\sigma _3 ^2 \Big)\label{EH}
 \eeq
with constant $B > 0, r \geq \sqrt[4]{B}, 0 \leq \psi \leq 2\pi$. The metrics are geodesically complete and asymptotically local Euclidean.
Topologically, the manifold is $$R_{\geq 0} \times SU(2)/Z_2 \cong R_{\geq 0} \times SO(3) \cong R_{\geq 0} \times P _3.$$

In \cite{Z}, metrics of Eguchi-Hanson type with the zero scalar curvature were constructed
 \beq
g =f^{-2} dr ^2+r ^2 \big( \sigma _1 ^2 +\sigma _2 ^2 +f^2 \sigma _3 ^2 \big) \label{EH-0}
 \eeq
with
 \beq
f=\sqrt{1-\frac{2A}{r ^2}-\frac{B}{r ^4}}, \label{f0}
 \eeq
where constant $B>0$, $r\geq \sqrt[4]{\frac{B}{n-1}}$ for $n \geq 2$, $0 \leq \psi \leq \frac{4\pi}{n}$ and
\beQ
A=-\frac{n-2}{2}\sqrt{\frac{B}{n-1}}.
\eeQ
The metrics are geodesically complete and asymptotically local Euclidean. Topologically, the manifold is $$R_{\geq 0} \times SU(2)/Z_n \cong R_{\geq 0} \times S^3/Z_n.$$

\mysection{Metrics of Eguchi-Hanson Type I}

In this section we constructed metrics of (\ref{EH-0}) with the constant scalar curvature, which we refer to metrics of Eguchi-Hanson type I.
The case of metrics of type I with the constant Ricci curvatures was studied extensively in \cite{P}. Let coframe of (\ref{EH-0}) be
\begin{equation*}
e^{1}=f^{-1}dr,\quad e^{2}=r \sigma_{1},\quad e^{3}=r \sigma_{2},\quad e^{4}=rf \sigma_{3}.
\end{equation*}
The connection 1-forms are defined by $$d e^{i}=-{\omega^{i}}_{j}\wedge e^{j}$$ and the curvature tensors are $${R^{i}} _{j}=d{\omega^{i}}_{j}+{\omega^{i}}_{k}\wedge{\omega^{k}}_{j}.$$ The scalar curvature $-24B$ satisfies the ordinary differential equation
\begin{equation*}
R=- (f ^2)'' -\frac{7}{r} (f ^2)' -\frac{8}{r^2}(f ^2 -1) =-24B, \label{R1}
\end{equation*}
where $B$ is certain constant. The general solutions are
\beq
f=\sqrt{1+\frac{C}{r^2}+\frac{A}{r^4}+Br^2} \label{f1}
\eeq
where $A$, $C$ are constant.

\begin{thm}\label{thm1}
Let constant $B>0$, natural number $n \geq 3$. Given any constant $C \leq \frac{(n-2)^2}{12B}$, metrics (\ref{EH-0}) with $f$ given by (\ref{f1}) and with ranges
\beq
r\geq r_0, \quad 0 \leq \psi \leq \frac{4\pi}{n} \label{range}
\eeq
are geodesically complete whose scalar curvature is $-24B$, where
\beq
r_0=\sqrt{\frac{n-2+\sqrt{(n-2)^2 -12BC}}{6B}}> 0,  \label{r0}
\eeq
\beq
A=\frac{1-2n +\sqrt{(n-2)^2 -12 BC}}{3}\,r_0 ^4. \label{A}
\eeq
Topologically, the manifold is $$R_{\geq 0} \times SU(2)/Z_n \cong R_{\geq 0} \times S^3/Z_n.$$
\end{thm}
\pf Let $\bar r$ be the largest positive root of $f$. Define
\beQ
\rho=r f \geq 0, \quad \rho(\bar r)=0.
\eeQ
By changing coordinates, metrics (\ref{EH-0}) with $f$ given by (\ref{f1}) become
\beq
\begin{aligned}
g =& \Big(f^2 +\frac{r}{2}(f^2)'\Big)^{-2} d\rho ^2 +\rho^{2}\sigma _3 ^2 +r ^2 \big( \sigma _1 ^2 +\sigma _2 ^2 \big)\\
  =& h^{-2} d\rho ^2 +\frac{\rho^2}{4} d\psi ^2 +\rho^2 \Big(\frac{\cos \theta}{2} d \phi d\psi +\frac{\cos^2 \theta}{4} d\phi ^2\Big)
  +r ^2 \big( \sigma _1 ^2 +\sigma _2 ^2 \big) \label{EH-1}
\end{aligned}
\eeq
where $h=1-\frac{A}{r^4}+2Br^2$. Now we find $\bar{r}$, $A$ such that
\beQ
\begin{aligned}
1+\frac{C}{\bar{r} ^2}+\frac{A}{\bar{r} ^4}+B \bar{r} ^2 =0,\quad
1-\frac{A}{\bar{r} ^4}+2B \bar{r} ^2=n.
\end{aligned}
\eeQ
We obtain that $\bar{r}=r_0$ satisfies (\ref{r0}) and $A$ satisfies (\ref{A}). Then metrics (\ref{EH-1}) become
\beQ
\begin{aligned}
g=& \frac{1}{n^2}\Big(d\rho ^2 +\rho^2 \big(d \frac{n\psi}{2} \big)^2 \Big)+\rho^2 \Big(\frac{\cos \theta}{2} d \phi d\psi +\frac{\cos^2 \theta}{4} d\phi ^2\Big)+r ^2 \big( \sigma _1 ^2 +\sigma _2 ^2 \big).
\end{aligned}
\eeQ
Thus the singularity at $r=r_0$ or $\rho =0$ can be removed under (\ref{range}) and the theorem holds.\qed

\begin{rmk}
The metrics in Theorem \ref{thm1} have natural K\"{a}hler structures \cite{L2}. But they are not ALH.
\end{rmk}

The curvature tensors of the metrics in Theorem \ref{thm1} are
\begin{equation*}
\begin{aligned}
&{R^{2}}_{1}={R^{3}}_{4}=\Big(-B+\frac{C}{r^4}+\frac{2A}{r^6}\Big)(e^2\wedge e^1+e^3\wedge e^4),\\
&{R^{3}}_{1}={R^{4}}_{2}=\Big(-B+\frac{C}{r^4}+\frac{2A}{r^6}\Big)(e^3\wedge e^1+e^4\wedge e^2),\\
&{R^{2}}_{3}=\Big(-4B-\frac{4C}{r^4}-\frac{4A}{r^6}\Big)e^2\wedge e^3+2\Big(-B+\frac{4C}{r^4}+\frac{2A}{r^6}\Big)e^4\wedge e^1,\\
&{R^{4}}_{1}=\Big(-4B-\frac{4A}{r^6}\Big)e^4\wedge e^1+2\Big(-B+\frac{4C}{r^4}+\frac{2A}{r^6}\Big)e^3\wedge e^2.
\end{aligned}
\end{equation*}
Thus the sectional curvatures are
\beQ
\begin{aligned}
K_{21}&=K_{31}=K_{34}=K_{42}=-B+\frac{C}{r^4}+\frac{2A}{r^6},\\
K_{23}&=-4B-\frac{4C}{r^4}-\frac{4A}{r^6},\\
K_{41}&=-4B-\frac{4A}{r^6}.
\end{aligned}
\eeQ
The Ricci curvatures are
\beQ
R_{11}=R_{44}=-6B+\frac{2C}{r^4},\quad R_{22}=R_{33}=-6B-\frac{2C}{r^4}.
\eeQ
If $C=0$, then the metrics are Einstein which were studied in \cite{P}.

In the following we show that the metrics constructed in Theorem \ref{thm1} can not develop into 5-dimensional spherically symmetric, static vacuum Einstein spacetimes with the negative cosmological constant
\beq
\tilde{R}_{ij}=\frac{2}{3}\Lambda \tilde{g} _{ij}, \quad 0\leq i,j \leq 4. \label{5D}
\eeq

\begin{thm}
There does not exist smooth function $v(r)$ such that the following metrics
\beq
\tilde{g}=-v^2 dt^2+f^{-2}dr^2+r^2(\sigma_{1}^2+\sigma_{2}^2+f^2\sigma_{3}^2)   \label{P5}
\eeq
satisfies (\ref{5D}).
\end{thm}
\pf Denote the coframe of $(\ref{P5})$
$$e^{0}=v dt,\quad e^{1}=f^{-1} dr,\quad e^{2}=r \sigma_{1},\quad e^{3}=r \sigma_{2},\quad e^{4}=rf\sigma_{3}.$$
It is easy to derive the connection 1-forms of $\tilde{g}$
\begin{equation*}
\begin{aligned}
{\tilde{\omega} ^{0}}_{1}&=\frac{v' f}{v}e^0,\,\,\, {\tilde{\omega} ^{0}}_{2}={\tilde{\omega} ^{0}}_{3}={\tilde{\omega} ^{0}}_{4}=0,\\
{\tilde{\omega} ^{2}}_{1}&=\frac{f}{r}e^2,\quad\, {\tilde{\omega} ^{3}}_{1}=\frac{f}{r}e^3,\quad\, {\tilde{\omega} ^{4}}_{1}=\Big(\frac{f}{r}+f'\Big)e^4,\\
{\tilde{\omega} ^{3}}_{4}&=\frac{f}{r}e^2,\quad\, {\tilde{\omega} ^{4}}_{2}=\frac{f}{r}e^3,\quad\, {\tilde{\omega} ^{2}}_{3}=\Big(\frac{2}{rf}-\frac{f}{r}\Big)e^4.
\end{aligned}
\end{equation*}
The curvature tensors of $\tilde{g}$ are
\begin{equation*}
\begin{aligned}
&{\tilde{R}^{0}}_{1}=-\Big(\big(1+\frac{C}{r^2}+\frac{A}{r^4}+Br^2\big)\frac{v''}{v}
+\big(\frac{-C}{r^3}+\frac{-2A}{r^5}+Br\big)\frac{v'}{v}\Big)e^0\wedge e^1,\\
&{\tilde{R}^{0}}_{2}=-\Big(\frac{1}{r}+\frac{C}{r^3}+\frac{A}{r^5}+Br\Big)\frac{v'}{v}e^0\wedge e^2,\\
&{\tilde{R}^{0}}_{3}=-\Big(\frac{1}{r}+\frac{C}{r^3}+\frac{A}{r^5}+Br\Big)\frac{v'}{v}e^0\wedge e^3,\\
&{\tilde{R}^{0}}_{4}=-\Big(\frac{1}{r}+\frac{-A}{r^5}+2Br\Big) \frac{v'}{v}e^0\wedge e^4,\\
&{\tilde{R}^{2}}_{1}={\tilde{R}^{3}}_{4}=\Big(-B+\frac{C}{r^4}+\frac{2A}{r^6}\Big)(e^2\wedge e^1+e^3\wedge e^4),\\
&{\tilde{R}^{3}}_{1}={\tilde{R}^{4}}_{2}=\Big(-B+\frac{C}{r^4}+\frac{2A}{r^6}\Big)(e^3\wedge e^1+e^4\wedge e^2),\\
&{\tilde{R}^{2}}_{3}=\Big(-4B-\frac{4C}{r^4}-\frac{4A}{r^6}\Big)e^2\wedge e^3+2\Big(-B+\frac{4C}{r^4}+\frac{2A}{r^6}\Big)e^4\wedge e^1,\\
&{\tilde{R}^{4}}_{1}=\Big(-4B-\frac{4A}{r^6}\Big)e^4\wedge e^1+2\Big(-B+\frac{4C}{r^4}+\frac{2A}{r^6}\Big)e^3\wedge e^2.
\end{aligned}
\end{equation*}
The Einstein equations (\ref{P5}) give
\begin{equation*}
\begin{aligned}
\tilde{R}_{00}&=-\big({\tilde{R} ^{0}}_{101}+{\tilde{R}^{0}}_{202}+{\tilde{R}^{0}}_{303}+{\tilde{R}^{0}}_{404}\big)=-\frac{2}{3}\Lambda,\\
\tilde{R}_{11}&=-6B+\frac{2C}{r^4}+{\tilde{R}^{0}}_{101}=\frac{2}{3}\Lambda,\\
\tilde{R}_{22}&=-6B-\frac{2C}{r^4}+{\tilde{R}^{0}}_{202}=\frac{2}{3}\Lambda,\\
\tilde{R}_{33}&=-6B-\frac{2C}{r^4}+{\tilde{R}^{0}}_{303}=\frac{2}{3}\Lambda,\\
\tilde{R}_{44}&=-6B+\frac{2C}{r^4}+{\tilde{R}^{0}}_{404}=\frac{2}{3}\Lambda.
\end{aligned}
\end{equation*}
From $\tilde{R}_{11}=\tilde{R}_{44}$, we obtain
\beQ
\Big(1+\frac{C}{r^2}+\frac{A}{r^4}+Br^2\Big)v''=\Big(\frac{1}{r}+\frac{C}{r^3}+\frac{A}{r^5}+Br\Big)v'.
\eeQ
Therefore
\beQ
v=\frac{C_1}{2}r^2+C_2.
\eeQ
From $\tilde{R}_{22}=\tilde{R}_{44}$, we obtain
\beQ
\Big(-\frac{C}{r^3}-\frac{2A}{r^5}+B r\Big)v'=\frac{4C}{r^4}v.
\eeQ
Substituting $v$ into the above equation, we obtain $C=C_1=0$ and $C_2$ is constant. Thus $v'=0$ which implies
\beQ
{\tilde{R} ^0}_{i0i}=0 \Longrightarrow \tilde{R} _{00}=0.
\eeQ
This gives contradiction. Therefore the proof of theorem is complete. \qed

\mysection{Metrics of Eguchi-Hanson Type II}\ls

In this section we study the following metrics with the negative constant scalar curvature. We refer them to metrics of Eguchi-Hanson type II.
\beq
g=\frac{dr^2}{(1+Br^2)f^2}+r^2\big(\sigma_{1}^2+\sigma_{2}^2+f^2\sigma_{3}^2 \big). \label{EH-2}
\eeq

Let $B>0$. When $f=1$, they are standard hyperbolic metrics. When
\beQ
f=\sqrt{1+\frac{A}{r^4}}, \quad A=-\frac{(n^2-4)^2}{16 B^2}, \quad n \geq 3
\eeQ
the metrics are geodesically complete and asymptotically local hyperbolic on
\beq
r\geq \sqrt{\frac{n^2 -4}{4B}}, \quad 0 \leq \psi \leq \frac{4\pi}{n}.
\eeq
They yield to the 5-dimensional Eguchi-Hanson-AdS spacetimes \cite{CM}
\beq
\tilde{g}=-(1+Br^2) dt^2+\frac{dr^2}{(1+Br^2)(1+\frac{A}{r^4})}+r^2\Big(\sigma_{1}^2+\sigma_{2}^2+(1+\frac{A}{r^4})\sigma_{3}^2\Big) \label{5D2}
\eeq
which satisfy the vacuum Einstein field equations (\ref{5D}). Furthermore, Dold \cite{D} calculated their Hawking mass at infinity
\beQ
m_H =-\frac{5}{6} \frac{(n^2-4)^2 }{16 B}<0.
\eeQ
He also studied the 5-dimensional Einstein vacuum equations and showed that the future maximal development does not contain horizon for small perturbation of Eguchi-Hanson type II initial data sets with negative Hawking mass at infinity.

Now we construct metrics with the negative constant scalar curvature. Let coframe of (\ref{EH-2}) be
\begin{equation*}
\begin{aligned}
e^{1}=(\sqrt{1+Br^2}f)^{-1}dr,\quad e^{2}=r \sigma_{1},\quad e^{3}=r \sigma_{2},\quad e^{4}=rf \sigma_{3}.
\end{aligned}
\end{equation*}
The connection 1-form $\left\{  {\omega^{i}}_{j}\right\}$ of (\ref{EH-2}) are
\begin{equation*}
\begin{aligned}
{\omega^{2}}_{1}&=\frac{\sqrt{1+Br^2}f}{r}e^2,\quad\ {\omega^{3}}_{1}=\frac{\sqrt{1+Br^2}f}{r}e^3,\\ {\omega^{4}}_{1}&=\Big(\frac{\sqrt{1+Br^2}f}{r}+\sqrt{1+Br^2}f'\Big)e^4,\\
{\omega^{3}}_{4}&=\frac{f}{r}e^2,\quad
{\omega^{4}}_{2}=\frac{f}{r}e^3, \quad
{\omega^{2}}_{3}=\Big(\frac{2}{rf}-\frac{f}{r}\Big)e^4.
\end{aligned}
\end{equation*}
The curvature tensors are
\begin{equation*}
\begin{aligned}
{R^{2}}_{1}=&-\Big(Bf^2+\frac{1+Br^2}{r}ff'\Big)e^2\wedge e^1+\frac{\sqrt{1+Br^2}ff'}{r}e^4\wedge e^3,\\
{R^{3}}_{1}=&-\Big(Bf^2+\frac{1+Br^2}{r}ff'\Big)e^3\wedge e^1+\frac{\sqrt{1+Br^2}ff'}{r}e^2\wedge e^4,\\
{R^{3}}_{4}=&-\Big(Bf^2+\frac{1+Br^2}{r}ff'\Big)e^3\wedge e^4+\frac{\sqrt{1+Br^2}ff'}{r}e^1\wedge e^2,\\
{R^{4}}_{2}=&-\Big(Bf^2+\frac{1+Br^2}{r}ff'\Big)e^4\wedge e^2+\frac{\sqrt{1+Br^2}ff'}{r}e^1\wedge e^3,\\
{R^{2}}_{3}=&-\Big(Bf^2+\frac{4f^2-4}{r^2}\Big)e^2\wedge e^3+\frac{2\sqrt{1+Br^2}ff'}{r}e^4\wedge e^1,\\
{R^{4}}_{1}=&-\Big(Bf^2+\frac{1+Br^2}{2}(f^2)''+\frac{3+4B r^2}{r}ff'\Big)e^4\wedge e^1 +\frac{2\sqrt{1+Br^2}ff'}{r}e^2\wedge e^3.
\end{aligned}
\end{equation*}
The constant scalar curvature $-12B$ satisfies the following equation
\beQ
R=-\Big((1+Br^2)(f^2)''+(\frac{7}{r}+8Br)(f^2)'+(\frac{8}{r^2}+12B)f^2-\frac{8}{r^2}\Big)=-12B.
\eeQ
The solutions are
\begin{equation}
f=\sqrt{1+\frac{\sqrt{1+Br^2}C}{r^4}+\frac{A}{r^4}}. \label{f2}
\end{equation}
where $A$, $C$ are constant.

Let $r_0$ be the largest positive root of $f$ given by (\ref{f2})
\beq
1+\frac{\sqrt{1+B r_0 ^2}C}{r _0 ^4}+\frac{A}{r _0 ^4}=0. \label{r02}
\eeq
It always exists, for example, if $A<0$. Choose $A$, $C$ such that
\begin{equation}
\sqrt{1+B r_0 ^{2} } \Big(1-\frac{A}{r_0 ^{4}}\Big) -\frac{C }{r_0 ^{4}} -\frac{B C}{2 r_0 ^{2} } = n. \label{n2}
\end{equation}
Then $r_0 ^2$ satisfies the following equation
\begin{equation}
\big(r_0 ^2 \big) ^{3}+\frac{4-n^2}{4B} \big(r_0 ^2 \big) ^{2}+\frac{n C}{4} \big(r_0 ^2 \big)-\frac{B C^{2}}{16}=0.  \label{t02}
\end{equation}

\begin{thm}\label{thm3}
Let constant $B>0$, natural number $n \geq 3$. If constant $C$ satisfies
\beq
\begin{aligned}
& C \leq \frac{n^3 - 36 n + (n^2 + 12) \sqrt{n^2 +12}}{27 B^2}\\
\mbox{or}\quad & C >\frac{(n^4 - 4)(3n + \sqrt{3 n^2 + 24})}{18 B^2}, \label{C2}
\end{aligned}
\eeq
then there exists $r_0>0$ such that metrics (\ref{EH-2}) with $f$ given by (\ref{f2}) and with ranges
\beQ
r\geq r_0, \quad 0 \leq \psi \leq \frac{4\pi}{n}
\eeQ
are geodesically complete and ALH whose scalar curvature is $-12B$, and
\beq
A=-r_0 ^4 -\sqrt{1+B r_0 ^2} C. \label{A2}
\eeq
Topologically, the manifold is
\beQ
R_{\geq 0} \times SU(2)/Z_n \cong R_{\geq 0} \times S^3/Z_n.
\eeQ
\end{thm}
\pf Define
\beQ
\rho=r f \geq 0, \quad \rho(r_0)=0.
\eeQ
By changing coordinates, metrics (\ref{EH-2}) with $f$ given by (\ref{f2}) become
\begin{equation*}
\begin{aligned}
g&=\frac{ d\rho ^2 }{(1+Br^2)\big(f^2 +\frac{r}{2}(f^2)' \big)^2}+\rho ^2 \sigma _3 ^2+r ^2 \big( \sigma _1 ^2 +\sigma _2 ^2 \big)\\
 &=h^{-2} d\rho ^2 +\rho ^2 \sigma _3 ^2+r ^2 \big(\sigma _1 ^2 +\sigma _2 ^2 \big),
\end{aligned}
\end{equation*}
where $$h(r)=\sqrt{1+B r^{2} } \Big(1-\frac{A}{r^{4}}\Big) -\frac{C}{r^{4}} -\frac{BC}{2 r^{2} }.$$ Denote
\beQ
\begin{aligned}
C_1&=\frac{n^3 - 36 n - (n^2 + 12) \sqrt{n^2 +12}}{27 B^2},\\
C_2&=\frac{n^3 - 36 n + (n^2 + 12) \sqrt{n^2 +12}}{27 B^2},\\
C_3&=\frac{(n^4 - 4)(3n - \sqrt{3 n^2 + 24})}{18 B^2},\\
C_4&=\frac{(n^4 - 4)(3n + \sqrt{3 n^2 + 24})}{18 B^2}.
\end{aligned}
\eeQ
It is obvious that $C_1<C_3<C_2<C_4$. Let
\beQ
\begin{aligned}
p&=\frac{-(n ^2 -4)^2 + 12 n B ^2 C}{48 B^2},\\
q&=\frac{-(n ^2 -4)^3 + 18 n (n^2 -4) B ^2 C - 54 B ^4 C ^2}{864 B ^3}.
\end{aligned}
\eeQ
A direct computation shows that
$$\Delta=\frac{p^3}{27}+\frac{q^2}{4}=\frac{27 B^4 C^4 -2 B^2 C^3 n (n ^2 -36 ) - 4 (n^2 - 4)^2 C^2 }{27648 B^2}.$$
We obtain
\begin{itemize}
\item[(i)] $C<C_1$ or $C>C_2$ $\Longrightarrow \Delta > 0$,\\
\item[(ii)] $C_1 \leq C \leq C_2$ $\Longrightarrow \Delta \leq 0$,\\
\item[(iii)] $C \leq C_3$ or $C\geq C_4$ $\Longrightarrow q \leq 0$.
\end{itemize}
Therefore, for any given $C$ satisfying (\ref{C2}), we solve (\ref{t02})
\begin{itemize}
\item[(i)] $C<C_1$ or $C>C_4$ $\Longrightarrow \Delta > 0, q<0$,
\beQ
r_0=\sqrt{ \sqrt[3]{-\frac{q}{2}+\sqrt{\frac{p^3}{27}+\frac{q^2}{4}}}+\sqrt[3]{-\frac{q}{2}-\sqrt{\frac{p^3}{27}+\frac{q^2}{4}}}+\frac{n^2 - 4}{12 B} },
\eeQ
\item[(ii)] $C_1 \leq C \leq C_2$ $\Longrightarrow \Delta \leq 0, p \leq 0$,
\beQ
r_0=\sqrt{ \frac{2 \sqrt{-3 p}}{3} \cos \frac{\alpha}{3} +\frac{n^2 - 4}{12 B} }
\eeQ
where $\alpha=\arccos \Big(\frac{-3 q \sqrt{-3 p}}{2 p^{2}}\Big)$, $0 \leq \alpha \leq \pi$.
\end{itemize}
With $A$ given by (\ref{A2}), the metrics become
\beQ
\begin{aligned}
g=& \frac{1}{n^2}\Big(d\rho ^2 +\rho^2 \big(d \frac{n\psi}{2} \big)^2 \Big)+\rho^2 \Big(\frac{\cos \theta}{2} d \phi d\psi +\frac{\cos^2 \theta}{4} d\phi ^2\Big)+r ^2 \big( \sigma _1 ^2 +\sigma _2 ^2 \big).
\end{aligned}
\eeQ
Therefore the theorem follows. \qed


The Ricci curvatures are
\begin{equation*}
\begin{aligned}
R_{11}=R_{44}=-3B-\frac{AB}{r^4},\quad R_{22}=R_{33}=-3B+\frac{AB}{r^4}.
\end{aligned}
\end{equation*}

\begin{coro}
Let $A=0$. The metrics constructed in Theorem \ref{thm3}
\beq
g=\frac{dr^2}{(1+Br^2)\big(1+\frac{\sqrt{1+Br^2} C}{r^4}\big)}+r^2\Big(\sigma_{1}^2+\sigma_{2}^2+\big(1+\frac{\sqrt{1+Br^2} C}{r^4}\big)
\sigma_{3}^2\Big)
\eeq
have constant Ricci curvature $-3B$.
\end{coro}

In the following we show that the metrics constructed in Theorem \ref{thm3} can not develop into 5-dimensional spherically symmetric, static vacuum Einstein spacetimes with the negative cosmological constant except (\ref{5D2}).

\begin{thm}
Let $f$ be given by (\ref{f2}). Suppose smooth function $v(r)$ such that the following metrics
\beq
\tilde{g}=-v^2 dt^2+\frac{dr^2}{(1+Br^2)f^2}+r^2\big(\sigma_{1}^2+\sigma_{2}^2+f^2\sigma_{3}^2\big)   \label{P5-2}
\eeq
satisfy (\ref{5D}). Then
\beQ
v=\sqrt{1+Br^2}, \quad C=0.
\eeQ
\end{thm}
\pf Denote $h=\sqrt{(1+B r^2)}f$. Then the curvature tensors are
\begin{equation*}
\begin{aligned}
{\tilde{R}^{0}}_{101}&=-\frac{v''h^2+v'h'h}{v},\\
{\tilde{R}^{0}}_{202}&={\tilde{R}^{0}}_{303}=-\frac{v'h^2}{vr},\\
{\tilde{R}^{0}}_{404}&=-(\frac{v'h^2}{vr}+\frac{v'f'h^2}{vf}),\\
{\tilde{R}^{2}}_{121}&={\tilde{R}^{3}}_{131}=-\frac{h'h}{r},\\
{\tilde{R}^{4}}_{141}&=-(\frac{h'h}{r}+\frac{f''h^2+f'h'h}{f}+\frac{2f'h^2}{rf}),\\
{\tilde{R}^{3}}_{434}&={\tilde{R}^{4}}_{242}=-(\frac{h^2-f^2}{r^2}+\frac{f'h^2}{rf}),\\
{\tilde{R}^{2}}_{323}&=-\frac{h^2+3f^2-4}{r^2}.
\end{aligned}
\end{equation*}
The Ricci curvatures are
	\begin{equation*}
	\begin{aligned}
	\tilde{R}_{00}&=-({\tilde{R}^{0}}_{101}+{\tilde{R}^{0}}_{202}+{\tilde{R}^{0}}_{303}+{\tilde{R}^{0}}_{404}),\\
	\tilde{R}_{11}&=-3B-\frac{AB}{r^4}+{\tilde{R}^{0}}_{101},\\
    \tilde{R}_{22}&=-3B+\frac{AB}{r^4}+{\tilde{R}^{0}}_{202},\\
	\tilde{R}_{33}&=-3B+\frac{AB}{r^4}+{\tilde{R}^{0}}_{303},\\
    \tilde{R}_{44}&=-3B-\frac{AB}{r^4}+{\tilde{R}^{0}}_{404}.
	\end{aligned}
	\end{equation*}
From $\tilde{R}_{11}=\tilde{R}_{44}$, we get
\begin{equation}
h^2v''-\Big(\frac{h^2}{r}+\frac{f'h^2}{f}-h'h\Big)v'=0.\label{v2}
\end{equation}
The general solutions of the (\ref{v2}) are
\beQ
v=\sqrt{1+B r^2}C_1+C_2.
\eeQ
From $\tilde{R}_{22}=\tilde{R}_{44}$, we get
\beQ
C(4+3Br^2)C_1-4AC_2=0.
\eeQ
So we obtain $A=0$, $C_1=0$ or $C=0$, $C_2=0$. In the first case, we obtain
\beQ
R_{00}=0.
\eeQ
This gives a contradiction. Hence, by changing $t$ to $C_1 t$, we can normalize $C_1=1$. Therefore
\beQ
v=\sqrt{1+B r^2}.
\eeQ
Thus we obtain (\ref{5D2}). \qed

\mysection{Energy}\ls

In this section we compute the total energy of ALH metrics of Eguchi-Hanson type II and show that it could be negative. For spacetimes with negative cosmological constant, the definition of total energy and other conserved quantities can be found in, e.g., \cite{HT, CN, W, CH, CMT, WXZ, WX} for four and higher dimensional spacetimes. We follow the method of \cite{HT} to provide the total energy related to metrics (\ref{EH-2}).

The 5-dimensional AdS spacetime with negative cosmological constant $-6B$ ($B>0$) is the hyperboloid
\begin{equation*}
 \eta_{\alpha \beta} y^{\alpha}y^{\beta}=-B^{-1}
\end{equation*}
in $R^{4,2}$ equipped with the metric
\begin{equation}
 ds^2=-\big(dy^{0}\big)^2+\sum_{i=1}^{4}\big(dy^{i}\big)^2-\big(dy^{5}\big)^2. \label{R42}
\end{equation}
Under coordinate transformations
\begin{equation*}
\begin{aligned}
y^0&=\frac{\cos(\sqrt{B}t)}{\sqrt{B}}\sqrt{1+B r^{2}},\quad y^1=r\cos\frac{\theta}{2}\cos\frac{\psi+\phi}{2},\\
y^2&=r\cos\frac{\theta}{2}\sin\frac{\psi+\phi}{2},\qquad \,\,\,y^3=r\sin\frac{\theta}{2}\cos\frac{\psi-\phi}{2},\\
y^4&=r\sin\frac{\theta}{2}\sin\frac{\psi-\phi}{2},\qquad \,\,\,y^5=\frac{\sin(\sqrt{B} t)}{\sqrt{B}}\sqrt{1+B r^{2}}.
\end{aligned}
\end{equation*}
The metric (\ref{R42}) reduces the 5-dimensional AdS metric
\begin{equation}
\tilde{g}_0 =-(1+B r^2)dt^2 +\frac{1}{1+B r^2}dr^2 +r^2\big(\sigma_{1}^2+\sigma_{2}^2+\sigma_{3}^2\big). \label{AdS5}
\end{equation}

Fifteen Killing vector fields
\begin{equation*}
U_{\alpha\beta}=y_{\alpha}\frac{\partial}{\partial y^{\beta}}-y_{\beta}\frac{\partial}{\partial y^{\alpha}}, \quad y_{\alpha}=\eta _{\alpha \beta} y^\beta
\end{equation*}
generate rotations for $R^{4,2}$ and reduce the symmetric structures of (\ref{R42}) restricted to the hyperboloid. The vector vector $U_{50}$ will generate the total energy.

Denote the hyperbolic metric
\begin{equation}
\breve{g}=\frac{1}{1+B r^2}dr^2 +r^2\big(\sigma_{1}^2+\sigma_{2}^2+\sigma_{3}^2\big) \label{H4}
\end{equation}
with sectional curvature $-B$. Let the coframe \{$\breve{e}^{\alpha}$\} of metric (\ref{AdS5}) be
\begin{equation*}
\begin{aligned}
\breve{e}^{0}=\sqrt{1+Br^2}dt,\quad \breve{e}^{1}=\frac{dr}{\sqrt{1+Br^2}},\quad \breve{e}^{2}=r\sigma_{1},\quad
\breve{e}^{3}=r\sigma_{2},\quad \breve{e}^{4}=r\sigma_{3}.
\end{aligned}
\end{equation*}
Let \{$\breve{e}_{\alpha}$\} be its dual frame. Then \{$\breve{e}^{i}$\} $(i=1,2,3,4)$ are also the coframe of the metric (\ref{H4}) and \{$\breve{e}_{i}$\} its dual frame.

The total energy of 4-dimensional ALH metrics with $S^3/Z_n$ topology on ends is defined as
\begin{equation}
E=\frac{1}{4\mbox{Vol}(S^3/Z_n)}\lim _{r \rightarrow \infty}\int _{S^3(r)/Z_n }\mathcal{E}
U_{50}^{(0)} \breve{\omega},     \label{E2}
\end{equation}
where $\breve{\nabla}$ is the Levi-Civita connection of $\breve{g}$, $g _{ij} =g\big(\breve{e}_i, \breve{e} _j\big)$, $a_{ij}=g_{ij}-\delta_{ij}$, $U_{50}=U_{50}^{\gamma} \breve{e}_{\gamma}$, $\breve{\omega}=\breve{e}^{2}\wedge\breve{e}^{3}\wedge\breve{e}^{4}$ and
\beQ
\mathcal{E}=\breve{\nabla}^{i} g_{1 i} - \breve{\nabla}_{1} \operatorname{tr}_{\breve g}(g)-\sqrt{B}\big(a_{11} -g_{11}\operatorname{tr}_{\breve g}(a)\big).
\eeQ

\begin{thm}
The metrics constructed in Theorem \ref{thm3} have the total energy $$E=A \sqrt{B}.$$ And $E<0$ if $A$ is chosen to be negative.
\end{thm}
\pf Direct computation shows that
\beQ
\begin{aligned}
\frac{\partial}{\partial t}&=\frac{\partial}{\partial y^0}\frac{\partial y^0}{\partial t}
+\frac{\partial}{\partial y^5}\frac{\partial y^5}{\partial t}
=-\sqrt{B} y^5 \frac{\partial}{\partial y^0}+\sqrt{B}y^0 \frac{\partial}{\partial y^5}=\sqrt{B} U_{50}.
\end{aligned}
\eeQ
Therefore
\beQ
U_{50}^{(0)}=\frac{\sqrt{1+Br^2}}{\sqrt{B}}.
\eeQ
The connection 1-forms of (\ref{H4}) are
\begin{equation*}
\begin{aligned}
{\breve \omega^{2}}_{1}&=\frac{\sqrt{1+Br^2}}{r}\breve e^2,\quad
{\breve \omega^{3}}_{1}=\frac{\sqrt{1+Br^2}}{r}\breve e^3,\quad
{\breve \omega^{2}}_{3}=\frac{1}{r}\breve e^4,\\
{\breve \omega^{4}}_{1}&=\frac{\sqrt{1+Br^2}}{r}\breve e^4,\quad
{\breve \omega^{4}}_{2}=\frac{1}{r}\breve e^3,\quad
{\breve \omega^{3}}_{4}=\frac{1}{r}\breve e^2.
\end{aligned}
\end{equation*}
The components of (\ref{EH-2}) are
\begin{equation*}
\begin{aligned}
g_{11}=f^{-2},\quad g_{22}=g_{33}=1,\quad g_{44}=f^2.
\end{aligned}
\end{equation*}
Therefore
\begin{equation*}
\begin{aligned}
\breve{\nabla}^{i} g_{1 i}&=\sum_{i=1}^{4}\left(\breve{e}_{i}\left(g_{i 1}\right)-g_{1 l} \breve{\omega}_{i}^{l}\left(\breve{e}_{i}\right)-g_{i l} \breve{\omega}_{1}^{l}\left(\breve{e}_{i}\right)\right)\\
&=\breve{e}_{1}\left(g_{1 1}\right)-g_{1 1} \breve{\omega}_{i}^{1}\left(\breve{e}_{i}\right)-g_{i i} \breve{\omega}_{1}^{i}\left(\breve{e}_{i}\right)\\
&=\sqrt{1+Br^2}\Big((f^{-2})'+\frac{3f^{-2}-f^2-2}{r}\Big),\\
\breve{\nabla}_{1} \operatorname{tr}_{\breve g}(g)&=\sqrt{1+Br^2}\Big((f^{-2})'+(f^2)' \Big).
\end{aligned}
\end{equation*}
Substituting (\ref{f2}) into the above formulas, we obtain
\beQ
\begin{aligned}
\frac{\sqrt{1+Br^2}}{\sqrt{B}}\Big(\breve{\nabla}^{i} g_{1 i}-\breve{\nabla}_{1} \operatorname{tr}_{\breve g}(g)\Big)&=-\frac{CB\sqrt{B}}{r\sqrt{1+Br^2}}+O(\frac{1}{r^4}),\\
\frac{\sqrt{1+Br^2}}{\sqrt{B}}\sqrt{B}\Big(a_{11}-g_{11}\operatorname{tr}_{\breve g}(a)\Big)
&=-\frac{BC}{r^2}-\frac{A\sqrt{1+Br^2}}{r^4}+O(\frac{1}{r^4}).
\end{aligned}
\eeQ
Thus
\beQ
\mathcal{E}U_{50}^{(0)}=\frac{A\sqrt{1+Br^2}}{r^4}+\frac{BC}{r^2 \big(1+Br^2 +\sqrt{(1+Br^2)Br^2}\big)}+O(\frac{1}{r^4}).
\eeQ
We obtain
\begin{equation*}
E=A\sqrt{B}.
\end{equation*}
From (\ref{C2}), (\ref{A2}), we find that $A$ is allowed to be negative. In this case $E<0$. \qed

\bigskip

{\footnotesize {\it Acknowledgement. The work of J. Chen is supported the special foundation for Guangxi Ba Gui Scholars. The work of X. Zhang is supported by Chinese NSF grants 11571345, 11731001, the special foundation for Guangxi Ba Gui Scholars, and HLM, NCMIS, CEMS, HCMS of Chinese Academy of Sciences.}

\end{document}